\documentclass[10pt]{article}

\begin{document}

\newtheorem{definition}[subsection]{Definition}
\newtheorem{theorem}[subsection]{Theorem}
\newtheorem{lemma}[subsection]{Lemma}
\newtheorem{cortje}[section]{Corollary}
\newtheorem{lemmatje}[section]{Lemma}
\newtheorem{proposition}[subsection]{Proposition}
\newtheorem{corollary}[subsection]{Corollary}
\newtheorem{proof}[subsection]{Proof}
\newtheorem{remark}[subsection]{Remark}
\newtheorem{observation}[subsection]{Observation}

\title{Large characteristic subgroups of surface groups not containing any simple loops}
\author{M. Pikaart }
\date{}
\maketitle

\begin{abstract}
We determine the largest (i.e. smallest index) characteristic
subgroup of surface groups not containing any simple loops.
\end{abstract}

\section{Introduction}
For any compact orientable surface $S$ we determine the smallest
characteristic non-geometric quotient of $\pi_1(S)$. Non-geometric
means that no non-trivial element that can be represented by a
simple closed curve is mapped to the identity and characteristic
means that its kernel is kept fixed under all automorphisms. We
write $\pi$, $H$ and $g$ for $\pi_1(S)$, $H_1(S, {\mathbf Z})$ and the genus of $S$
respectively.
We assume that $g$ is at least 2.
Consider the following characteristic subgroups of
$\pi$ : $\pi^{[1]}  :=  \pi$ and inductively $ \pi^{[k+1]}  :=
[\pi,\pi^{[k]}]$. We have a well known isomorphism
$\pi^{[2]}/\pi^{[3]} \rightarrow \wedge^2 H / \omega$, $[x,y]
\mapsto x\wedge y$, where $\omega $ is the intersection form on
$H$. We have an intersection product $\wedge^2 H
\rightarrow {\mathbf Z} $.
 Let $K$ be the kernel of the composition of the map $\pi^{[2]}
 \rightarrow \wedge^2 H / \omega$ with the intersection product $\wedge^2
 H/\omega \rightarrow {\mathbf Z}/g{\mathbf Z} $.
 Our result is
\begin{proposition}
If $g$ is odd, then the largest (i.e. smallest index)
characteristic non-geometric subgroup of $\pi$ is given by all
$g$th powers and $K$. If $g$ is even, then the largest
characteristic non-geometric subgroup of $\pi$ is given by all
$2g$th powers and $K$. The indices of these are $g^{2g+1}$ and
$(2g)^{2g}g$ respectively.
\end{proposition}

Recently, Livingston proved that for $g=2$ the smallest
non-geometric quotient of $\pi$ is a group of order $2^5$, and he
raised the question whether an easy generalisation of his result
(which yields a group of order $g^{2g+1}$), holds true for any
genus. The theorem above shows that his generalisation for odd
genus is in any case the smallest non-geometric characteristic
quotient.

\section{Notation and preliminary computations}

Before we
start, we pick once and for all a set of generators for $\pi$:
$x_1, \dots, x_{2g}$ with defining relation $\Pi_{i=1}^g
[x_{2i-1},x_{2i}]$. We call a pair of integers $(i,j)$ related if
and only if there exists an integer $h$ with $1 \leq h \leq g-1$
such that $(i,j)=(2h-1,2h)$.
 We write $[a,b]$ for $a^{-1}b^{-1}ab$, such
that $ab=ba[a,b]$.
Furthermore, if $\alpha$ and $\beta$ are elements of $\pi$, then
$\alpha \beta$ means: perform $\alpha$ first, then $\beta$.

Consider the following normal subgroups of $\pi$:
$$\begin{array}{rcl}
    \pi^{[1]} & := & \pi; \\
     \pi^{[k+1]} & := & [\pi,\pi^{[k]}];\\
    K& := & Ker(\pi^{[2]} \rightarrow {\mathbf Z}/g {\mathbf Z}); \\
     \pi^n & := & <x^n | x \in \pi>.\\

\end{array} $$

Notice that $K$ contains $\pi^{[3]}$. All these are clearly
characteristic subgroups. For $K$ this holds since $\pi^{[2]}$ is
characteristic and all involved maps are natural. The quotient we
are interested in, is $\pi/\pi^g K$ respectively $\pi/\pi^{2g}K$.
In order to prove the proposition, we first describe
$\pi/\pi^{[3]}$.

\begin{lemma}
$\pi/\pi^{[3]} \cong \{ \Pi_{i=1}^{2g} x_i^{n_i}\ \Pi_{1 \leq i <
j \leq 2g}^* [x_i,x_j]^{m_{i,j}} | n_i,m_{i,j} \in {\mathbf Z},
*\ meaning: (i,j) \neq (2g-1,2g) \},$

where multiplication on the right hand side is defined as follows:
$$ \begin{array}{c} (\Pi_{i=1}^{2g} x_i^{n_i} \Pi_{1 \leq i < j \leq
2g}^* [x_i,x_j]^{m_{i,j}}) \ ( \Pi_{i=1}^{2g} x_i^{k_i} \Pi_{1 \leq
i < j \leq 2g}^* [x_i,x_j]^{l_{i,j}})\\

\\ =
\Pi_{i=1}^{2g}x_i^{n_i+k_i} \Pi_{1 \leq i < j \leq 2g}^*
[x_i,x_j]^{m_{i,j}+l_{i,j}-k_in_j+\tilde{\delta}_{i+1,j}k_{2g-1}n_{2g}}.
\end{array} $$
Here $\tilde{\delta}_{a,b}=1$ if $a=b$ and $a$ is even, else
$\tilde{\delta}_{a,b}=0$.
\end{lemma}

\begin{proof}{\rm
We have the exact sequence $$1 \rightarrow \pi^{[2]}/\pi^{[3]}
\rightarrow \pi/ \pi^{[3]} \rightarrow \pi / \pi^{[2]} \rightarrow
1 ,$$ where the outer factors are finitely generated free abelian
groups (\cite[Th.5.12]{Labute}) of rank $(^{2g}_2)-1$ and $2g$
respectively. Explicitly, we have:
\[ \begin{array}{l}
 \pi / \pi^{[2]} = \langle x_i\rangle_{ab} ~,i=1, \dots, 2g, \\
 \pi^{[2]}/\pi^{[3]}= \langle
[x_i,x_j]\rangle_{ab}~, 1 \leq i < j \leq 2g, ~(i,j) \neq
(2g-1,2g).
\end{array} \]
Here the subscript ab means the free abelian group generated by
these elements,
Now recall that if $a \in \pi^{[k]}$ and $b \in \pi^{[l]}$ then
$ab=ba[a,b]$
 with $[a,b] \in \pi^{(k+l)}$, so that
\begin{eqnarray} ab\equiv ba \mbox{  mod } \pi^{[k+l]}. \end{eqnarray}
Furthermore we have the following identities, modulo $\pi^{[3]}$:
(\cite[Th.5.1]{Magnus})
\begin{eqnarray}
\ [a,b]  & = & [b,a]^{-1},\\ \ [a,bc] & = & [a,c]~[a,b],
\\ \ [ab,c] & = & [a,c]~[b,c].
\end{eqnarray}
It is clear now that any element of $\pi / \pi^{[3]}$ can be
written uniquely in the form $ \Pi_{i=1}^{2g} x_i^{n_i}  \Pi_{1 \leq i < j
\leq 2g}^* [x_i,x_j]^{m_{i,j}}$
where $ n_i,m_{i,j} \in {\mathbf Z}$ and the meaning of the * is: $ (i,j) \neq (2g-1,2g)$.
Note that the last
$(^{2g}_2)-1$ factors commute by (1). Before we start the
computation, notice that modulo $\pi^{[3]}$ we have
\begin{eqnarray}
[a^i,b^j]\equiv [a,b]^{ij},
\end{eqnarray}
as one proves easily by induction. The following identities hold
in the group $\pi / \pi^{[3]}$:
$$ \begin{array}{cl}
&(x_1^{n_1} \dots x_{2g}^{n_{2g}})
(x_1^{k_1} \dots x_{2g}^{k_{2g}})\\
=& x_1^{n_1+k_1} \dots x_{2g}^{n_{2g}}
x_2^{k_2} \dots x_{2g}^{k_{2g}}[x_1,x_2]^{-k_1n_2} \dots
[x_1,x_{2g}]^{-k_1n_{2g}}\\
=& \dots \\
=&\Pi x_i^{n_i+k_i}\Pi_{1 \leq i < j \leq 2g} [x_i,x_j]^{-k_in_j}\\
=& \Pi x_i^{n_i+k_i}\Pi_{1 \leq i < j \leq 2g}^*
[x_i,x_j]^{-k_in_j+\tilde{\delta}_{i+1,j}k_{2g-1}n_{2g}}
\end{array} $$
with $\tilde{\delta}_{i+1,j}k_{2g-1}n_{2g}$ as above. Combining
these proves the lemma. (Cf.\cite[Lemma 6.1]{PdJ})}
\end{proof}
 Here the term with
$\tilde{\delta}_{i+1,j}k_{2g-1}n_{2g}$ stems from the defining
relation for the group $\pi: \Pi_{i=1}^{2g} [x_{2i-1},x_{2i}]=1$,
which we used to get rid of $[x_{2g-1},x_{2g}]$ as a generator for
$\pi^{[2]}/ \pi^{[3]}$.

Clearly the subgroup $\pi^{[2]} / \pi^{[3]}$ is given
in terms of these generators by all expressions of the form
$\Pi_{1 \leq i < j \leq 2g}^* [x_i,x_j]^{m_{i,j}}$.

\begin{lemma}
$\pi/K \cong \{ \Pi_{i=1}^{2g} x_i^{n_i}\
[x_1,x_2]^m | n_i \in {\mathbf Z}, m \in {\mathbf Z}/ g{\mathbf Z} \},$
where multiplication on the right hand side is defined as follows:
$$ \begin{array}{c} (\Pi_{i=1}^{2g} x_i^{n_i} [x_1,x_2]^m)\ ( \Pi_{i=1}^{2g} x_i^{k_i}
[x_1,x_2]^l)\\

\\ =
\Pi_{i=1}^{2g}x_i^{n_i+k_i}
[x_1,x_2]^{m+l-\Sigma_{i=1}^{g-1}k_{2i-1}n_{2i}+(g-1)k_{2g-1}n_{2g}}.
\end{array} $$
\end{lemma}

\begin{proof}{\rm
Clearly, $K / \pi^{[3]}$ is generated by the elements $[x_i,x_j]$
for $i$ and $j$ not related, by the elements
$[x_{2k-1},x_{2k}][x_{2l-1},x_{2l}]^{-1}$ and by $\Pi_{i=1}^g [x_{2i-1},x_{2i}]$.
}
\end{proof}

\begin{proposition}\label{beschrijving}
If $g$ is odd, then
$$\pi/\pi^g K  \cong \{ \Pi_{i=1}^{2g} x_i^{n_i}[x_1,x_2]^m |
 n_i,m \in {\mathbf Z}/ g{\mathbf Z}, \},$$
and if $g$ is even, say $g=2h$, then
$$\pi/\pi^g K  \cong \{ \Pi_{i=1}^{2g} x_i^{n_i}[x_1,x_2]^m |
 n_i \in {\mathbf Z}/ g{\mathbf Z},m \in {\mathbf Z}/ h{\mathbf Z}
 \},$$
where multiplication on the right hand sides is defined as above.
\end{proposition}

\begin{proof}{\rm
Consider the short exact sequence
$$ 1 \rightarrow
 \pi^{[2]}/(\pi^{[2]} \cap \pi^g K ) \rightarrow \pi/ \pi^g K
\rightarrow \pi/ \pi^g \pi^{[2]} \rightarrow 1.$$

Clearly, $\pi / \pi^g \pi^{[2]}\cong H_1(S, {\mathbf Z}/g{\mathbf Z})$.
Furthermore, if $g$ is odd,
 $ \pi^{[2]}/(\pi^{[2]} \cap \pi^g K )  \cong
 ({\mathbf Z}/g{\mathbf Z})$, whereas in case $g$ is even, say $g=2h$, we have
 $ \pi^{[2]}/(\pi^{[2]} \cap \pi^g K )  \cong
 ({\mathbf Z}/h{\mathbf Z})$ . This follows directly from \cite[Lemma
 6.3]{PdJ}.}
\end{proof}

\begin{corollary}
The subgroup $\pi^g$ is generated by all $g$th powers of geometric
(i.e. representable by a simple closed curve) elements modulo $K$.
\end{corollary}

\begin{proof}{\rm
Clearly the elements $x_i^g$ and $[x_1,x_2]^g$ are $g$th powers
of geometric elements, settling the statement for odd $g$. For
even $g$, say $g=2h$, we have $(x_1x_2)^{2h} \cong
x_1^{2h}x_2^{2h}[x_1,x_2]^{(2h-1)h}$ modulo $\pi^{[3]}$, as one proves easily by
induction (cf. \cite[Lemma 6.3]{PdJ}). Therefore, $[x_1,x_2]^{-h}
\cong x_2^{-2h}x_1^{-2h}(x_1x_2)^{2h}[x_1,x_2]^{-2h}$, proving the
corollary.}
\end{proof}

We define the following simple loops on $S$:
    $\gamma_1=x_1$,
    $\gamma_{2h}=x_{2h}$, for $h=1, \dots, g$,
    $\gamma_{2h-1}=x_{2h-1}[x_{2h-3},x_{2h-2}]^{-1}x_{2h-3}^{-1}$, for $h=2,
    \dots,g$,and finally
    $\gamma_{2g+1}$ $ = [x_{2g-1},x_{2g}]^{-1}x_{2g-1}^{-1}.$

We write $\tau_i$ respectively
$\sigma_i$ for (the right handed) Dehn twist around $\gamma_i$
and $x_{2i-1}$ respectively. For later convenience we list the
action of these Dehn twists on the generators of $\pi$  and the
action modulo $\pi^{[3]}$ as above (if the
action of some Dehn twist on a generator is not given, it is the
trivial action): $$
\begin{array}{rcl}
    \tau_1(x_2) & = & x_1^{-1}x_2, \\
    \tau_{2h}(x_{2h-1}) & = & x_{2h}x_{2h-1} \cong
x_{2h-1}x_{2h}[x_{2h-1},x_{2h}]^{-1}, \\
    \tau_{2h-1}(x_{2h-2}) & = & x_{2h-2}\gamma_{2h-1} =
x_{2h-2}x_{2h-1}[x_{2h-3},x_{2h-2}]^{-1}x_{2h-3}^{-1}\\ & \cong  &
x_{2h-3}^{-1}x_{2h-2}x_{2h-1}[x_{2h-3},x_{2h-2}]^{-2}[x_{2h-2},x_{2h-1}]^{-1}\\
\tau_{2h-1}(x_{2h-1}) & = &
\gamma_{2h-1}^{-1}x_{2h-1}\gamma_{2h-1}=
x_{2h-1}[x_{2h-1},\gamma_{2h-1}] \\ & = & x_{2h-1}[x_{2h-1},
x_{2h-1}[x_{2h-3},x_{2h-2}]^{-1}x_{2h-3}^{-1}] \\ & \cong &
x_{2h-1}[x_{2h-3}, x_{2h-1}] \\ \tau_{2h-1}(x_{2h}) & = &
\gamma_{2h-1}^{-1}x_{2h} \\
  & = & x_{2h-3}[x_{2h-3},x_{2h-2}]x_{2h-1}^{-1} x_{2h}\\
& \cong &  x_{2h-3}x_{2h-1}^{-1}x_{2h}[x_{2h-3}, x_{2h-2}] \\
\tau_{2g+1}(x_{2g}) & = & x_{2g}\gamma_{2g+1}=
x_{2g}[x_{2g-1},x_{2g}]^{-1}x_{2g-1}^{-1} \\
 & = & x_{2g-1}^{-1}x_{2g}\\
\sigma_i(x_{2i}) & = & x_{2i-1}^{-1}x_{2i}\\
\end{array} $$

\section{The proofs}

\begin{proposition}
The quotient $\pi/ \pi^g K$ is characteristic and it is
non-geometric if and only if $ g$ is odd. If $g$ is even, the
quotient $\pi/\pi^{2g} K $ is characteristic and non-geometric.
\end{proposition}

\begin{proof}{\rm
$\pi^g K $ is characteristic, while $\pi^g$ and $K$ are.
 For $g$ odd, we have that $\pi^gK$ is non-geometric since it
 is exactly the group described by Livingston in \cite[Section
 5]{Liv}. For $g$ even, say $g=2h$, it follows from Proposition \ref{beschrijving}
 that $\Pi_{i=1}^h [ x_{2i-1},x_i]$ is contained in $\pi^g K$, so
 $\pi^g K$ is geometric. On the other hand,  for $g$ even, Livingston's group is a quotient of
 $\pi^{2g}K$, namely the one generated by the elements $x_i^g$
 for $i=1, \dots, 2g$. Thus $\pi^{2g} K$ is non-geometric.}
\end{proof}

\begin{theorem}
If $g$ is odd, the largest characteristic non-geometric quotient of
$\pi$ is $\pi/\pi^g K $ and if $g$ is even, the largest characteristic
non-geometric quotient of $\pi$ is $\pi/\pi^{2g} K $. The
indices of these groups are $g^{2g+1} $ respectively $(2g)^{2g} g=(2g)^{2g+1}/2$.

\end{theorem}

\begin{proof}{\rm
The indices follow directly from Proposition \ref{beschrijving}.
Let $M$ be any characteristic finite-index non-geometric subgroup
of $\pi$. Let $k$ be the smallest positive integer such that for
some simple closed not separating curve $\delta$ we have $\delta^k
\in M$. Since $M$ is characteristic and all simple closed not
separating curves can be mapped one onto the other, we have that
$k$ is the smallest positive integer with $\delta^k \in M$ for all
simple closed not separating curves $\delta$. Notice that $k \geq
3$. Namely, if $k=2$, then $[x_1,x_2]=x_1^{-2} (x_1x_2)^2 x_2^{-2}
\in M$, so $M$ is geometric, contradiction. We have $\pi
/M\pi^{[2]} \cong ({\mathbf Z}/k {\mathbf Z})^{2g}$ for some
positive integer $k$.

Let $P$ be the subgroup of $\pi^{[2]}$ generated by all
$[x_i,x_j]$ for $(i,j)$ not related (recall that $i<j$). By an
analogous argument we have that $(\pi^{[2]} \cap M) / (\pi^{[3]}
\cap M)P$ is a quotient of $({\mathbf Z}/m {\mathbf Z})^{g-1}$,
generated by the elements $[x_{2i-1},x_{2i}]$ for $i=1,\dots ,g$.
Again by the classification of surfaces and by the fact that $M$
is characteristic, we get a uniform power $l$ such that
$[x_{2i-1},x_{2i}]^t \in M$ if and only if $l$ divides $t$. We
have that $$[x_{2i-1}^k,x_{2i}] \cong [x_{2i-1},x_{2i}]^k \mbox{
modulo }\pi^{[3]}$$ and thus $l$ divides $k$ if $k$ odd and $2l$
divides $k$ if $k$ even, by Proposition \ref{beschrijving}.
Similarly, we obtain a uniform power $m$ such that for $(i,j)$ not
related, the elements $[x_i,x_j]^t$ are in $M$ if and only if $m$
divides $t$. Furthermore we have that $m$ divides $l$, since
$\tau_3([x_1,x_2]^l) \cong [x_1,x_2]^l[x_1,x_3]^l$ modulo
$\pi^{[3]}$.

Now suppose $[x_1,x_2] \not \equiv [x_3,x_4]$ modulo $M$. Since
$M$ is characteristic and all classes $[x_{2i-1},x_{2i}] $ can be
transformed one into the other by an automorphism of $\pi$, it
follows that  they all have different classes modulo $M$, thus
$[x_{2i-1},x_{2i}][x_{2j-1},x_{2j}]^{-1}$ is not contained in $M$.
Consider the short exact sequence $$ 1 \rightarrow M \pi^{[3]}/
\pi^{[3]} \rightarrow \pi/ \pi^{[3]} \rightarrow \pi/M \pi^{[3]}
\rightarrow 1.$$ We claim that $M \pi^{[3]}/ \pi^{[3]}$ does not
contain any element of the form $[x_i,x_j]$ for $(i,j)$ not
related (equivalently, $m>1$). Namely, suppose that there is a
pair $(i,j)$, with $(i,j)$ not related, such that $[x_i,x_j] \in
M$. Then $[x_i,x_j]$ is in $M$ for all $(i,j)$ not related, again
by the classification of surfaces and the fact that $M$ is
characteristic. We compute $\tau_{2h-1}([x_{2h-2},x_{2h}])$. We
have the following identities modulo $\pi^{[3]}$, where $2 \leq h
\leq g$: $$ \begin{array}{cl} & \tau_{2h-1}([x_{2h-2},x_{2h}]) \\
 \cong &
[x_{2h-3}^{-1}x_{2h-2}x_{2h-1},x_{2h-3}x_{2h-1}^{-1}x_{2h}] \\
 \cong &[x_{2h-3},x_{2h-1}][x_{2h-3},x_{2h}]^{-1}
[x_{2h-2},x_{2h-1}]^{-1}[x_{2h-2},x_{2h}][x_{2h-3},x_{2h-1}]^{-1}
\\
& [x_{2h-3},x_{2h-2}]^{-1}[x_{2h-1},x_{2h}]\\
  \cong & [x_{2h-3},x_{2h}]^{-1}
[x_{2h-2},x_{2h-1}]^{-1}[x_{2h-2},x_{2h}]
[x_{2h-3},x_{2h-2}]^{-1}[x_{2h-1},x_{2h}].
\end{array} $$
The product
$[x_{2h-3},x_{2h-2}]^{-1}[x_{2h-1},x_{2h}]$ is contained in $M$
since the first three commutators are in $M$. This leads to a contradiction.

Now we claim that the abelian $m$-torsion subgroup of $\pi/ M \pi^{[3]}$ generated by all
$[x_i,x_j]$ for $(i,j)$ not related, has rank $(^{2g}_2)-g=2g^2-2g$.
Namely, suppose there is an element $z=\sum_{(i,j) not rel.} n_{i,j}
[x_i,x_j] $ contained in $M$ with all $n_{i,j} \in \{ 0,\dots, m-1 \}$
(in additive notation).
We compute a number of elements of the form $\tau(z)-z$ to show
that all these $n_{i,j}$ are actually zero. We use repeatly that
$n_{i,j}[x_i,x_j]$ is in $M$ if and only if $m $ divides $n_{i,j}$
$$ \begin{array}{rllll}
z_{2g}  := & \tau_2(z)-z &=  \sum_{3 \leq j \leq 2g} n_{1,j} [x_1,x_j] & &\\
z_{2g-1}  := & \sigma_g(z_{2g})-z_{2g} &=   n_{1,2g}[x_1,x_{2g-1}] & \Rightarrow & n_{1,2g}=0 \\
z_{2g-2}  := & \tau_{2g}(z_{2g-1})-z_{2g-1} &=   n_{1,2g-1}[x_1,x_{2g}] & \Rightarrow & n_{1,2g-1}=0 \\
z_{2g-3}  := & \sigma_{g-1}(z_{2g-2})-z_{2g-2} &=   n_{1,2g-2}[x_1,x_{2g-3}] & \Rightarrow & n_{1,2g-2}=0 \\
z_{2g-4}  := & \tau_{2g-2}(z_{2g-3})-z_{2g-3} &=   n_{1,2g-3}[x_1,x_{2g-2}] & \Rightarrow & n_{1,2g-3}=0 \\
 & & \dots &  & \\
 & \tau_4(z_{..})-z_{..} & = n_{1,3} [x_1,x_3] & \Rightarrow & n_{1,3} =0 \\
&  \mbox {we get }z& = \sum_{2 \leq j \leq 2g} n_{1,j} [x_1,x_j] \\
z_{2g}  := & \tau_1(z)-z &=  \sum_{3 \leq j \leq 2g} n_{2,j} [x_2,x_j] & &\\
\end{array} $$
It is clear that continuing in this way we we show that all
coefficients $n_{i,j}$ are zero. Since there are $(^{2g}_2)-2g$
elements of this form, this proves the claim.

It remains to show that the index of this quotient is larger than
the index of $\pi^g K $ respectively $\pi^{2g} K $. We have that
the index of $M$ is at least $\# (\pi/M \pi^{[2]}) \# (M \cap
\pi^{[2]}) /(M \cap \pi^{[3]})$, which in turn is at least
 $k^{2g} m^{2g^2-2g}$.
If $m$ is even, this leaves $m=2$ as smallest possibility. Since $m$ divides $l$,
$l$ is also even and therefore $2l$ divides $k$. So the smallest
possibility we get is $4^{2g}2^{2g^2-2g}=(2^g)^{2g+2}$. This is
larger than $(2g)^{2g+1}/2$ for all $g$.

If $m$ is odd, the smallest possibility becomes
$3^{2g}3^{2g^2-2g}=(3^g)^{2g}$, which is again larger than
$g^{2g+1}$ for all $g$.

On the other hand, suppose $[x_1,x_2] \equiv [x_3,x_4]$ modulo
$M$. By the same argument, we get that all the classes
$[x_{2i-1},x_{2i}] $ are equal modulo $M$. Since $\Pi_{i=1}^h
[x_{2i-1},x_{2i}]  \notin M$ for all $h=1, \dots, g-1$, the
smallest possibility for $(\pi^{[2]} \cap M) / (\pi^{[3]} \cap
M)P$ is $({\mathbf Z}/g {\mathbf Z})$.
This implies that $g$ divides $k$ and $2g$ divides $k$ if $g$
even. Thus, $M$ is contained in $\pi^g K $ or $\pi^{2g} K$
if $g$ is even.}

\end{proof}


\end{document}